\documentclass{amsart}

\usepackage{amssymb,latexsym,amsmath}

\begin{document}

\def\cc{{\mathbb C}}
\def\r{{\mathbb R}}
\def\q{{\mathbb Q}}
\def\z{{\mathbb Z}}
\def\i{{\mathbb I}}

\def\p{{\mathcal P}}
\def\h{{\mathcal H}}
\def\alg{A_{aut}(X_n)}
\def\no{|\! |}

\newtheorem{theo}{Theorem}[section]
\newtheorem{prop}{Proposition}[section]
\newtheorem{coro}{Corollary}[section]
\newtheorem{lemm}{Lemma}[section]
\theoremstyle{definition}
\newtheorem{defi}{Definition}[section]
\theoremstyle{remark}
\newtheorem{exam}{Example}[section]
\newtheorem{rema}{Remark}[section]

\title{On the structure of quantum permutation groups}
\author{Teodor Banica}
\address{Laboratoire Emile Picard, Universit\'e Paul Sabatier, 
118 route de Narbonne, 31062 Toulouse, France}
\email{banica@picard.ups-tlse.fr}
\author{Sergiu Moroianu}
\address{Institutul de Matematic\u a al Academiei Rom\^ane, P.O. Box 1-764,
RO-014700 Bucure\c sti, Romania}
\email{moroianu@alum.mit.edu}
\thanks{Moroianu was partially supported by 
the Marie Curie grant MERG-CT-2004-006375
funded by the European Commission, and by a CERES contract (2004)}
\date{\today}
\begin{abstract}
The quantum permutation group of the set $X_n=\{ 1,\ldots ,n\}$
corresponds to the Hopf algebra $A_{aut}(X_n)$. This is an algebra constructed with generators
and relations, known to be isomorphic to $\cc
(S_n)$ for $n\leq 3$, and to be infinite dimensional for
$n\geq 4$. In this paper we find an explicit representation of the algebra $A_{aut}(X_n)$,
related to Clifford algebras. For $n=4$ the representation is faithful in 
the discrete quantum group sense.
\end{abstract}
\maketitle

\section*{Introduction}

A general theory of unital Hopf $\cc^*$-algebras is developed by Woronowicz in \cite{wo1}, \cite{wo2}, \cite{wo3}. The main results are the existence of the Haar functional, an analogue of Peter-Weyl theory and of Tannaka-Krein duality, and explicit formulae for the square of the antipode. As for examples, these include algebras of continuous functions on compact groups, $q$-deformations of them with $q>0$, and $\cc^*$-algebras of discrete groups.

Of particular interest is the algebra $\alg$ constructed by Wang in \cite{wa2}. This is the universal Hopf $\cc^*$-algebra coacting on the set $X_n=\{ 1,\ldots ,n\}$. In other words, the compact quantum group associated to it is a kind of analogue of the symmetric group $S_n$.

The algebra $\alg$ is constructed with generators and relations. There are $n^2$ generators, labeled $u_{ij}$ with $i,j=1,\ldots ,n$. The relations are those making $u$ a magic biunitary matrix. This means that all coefficients $u_{ij}$ are projections, and on each row and each column of $u$ these projections are mutually orthogonal, and sum up to $1$. The comultiplication is given by $\Delta (u_{ij})=\sum u_{ik}\otimes u_{kj}$ and the fundamental coaction is given by $\alpha (\delta_i)=\sum u_{ji}\otimes\delta_j$.

For $n=1,2,3$ the canonical quotient map $\alg\to\cc (S_n)$ is an isomorphism. For $n\geq 4$ the algebra $\alg$ is infinite dimensional, and just a few things are known about it. Its irreducible corepresentations are classified in \cite{aut}, with the conclusion that their fusion rules coincide with those for irreducible representations of $SO(3)$, independently of $n\geq 4$. In \cite{wa4} Wang proves that the compact quantum group associated to $\alg$ with $n\geq 4$ is simple. In \cite{aut} it is shown that the discrete quantum group associated to $\alg$ with $n\geq 5$ is not amenable. Various quotients of $\alg$, corresponding to quantum symmetry groups of polyhedra, colored graphs etc., are studied in \cite{gr} by using planar algebra techniques.

These results certainly bring some light on the structure of $\alg$. However,
for $n\geq 4$ this remains an abstract algebra, constructed with generators and
relations.
 
In this paper we find an explicit representation of $A_{aut}(X_n)$. The
construction works when $n$ is a power of $2$, and uses a magic biunitary matrix related to Clifford algebras. For $n=4$ the representation is inner
faithful, in the sense that the corresponding unitary
representation of the discrete quantum group associated to  $A_{aut}(X_4)$ is faithful.

As a conclusion, there might be a geometric
interpretation of Hopf algebras of type $A_{aut}(X_n)$. We should
mention here that for the algebra $A_{aut}(X)$ with $X$ finite
graph such an interpretation would be of real
help, for instance in computing fusion rules.

\subsection*{\bf Acknowledgments}
We are deeply grateful to Georges Skandalis for all the discussions and his
essential contribution to the results of this paper.

\section{Magic biunitary matrices}

Let $A$ be a unital $\cc^*$-algebra. That is, we have a unital algebra $A$ over the field of complex numbers $\cc$, with an antilinear antimultiplicative map $a\to a^*$ satisfying $a^{**}=a$, and with a Banach space norm satisfying $\no a^*a\no =\no a\no^2$.

A projection is an element $p\in A$ satisfying $p^2=p^*=p$. Two projections $p,q$ are said to be orthogonal when $pq=0$. A partition of unity in $A$ is a finite set of projections, which are mutually orthogonal, and sum up to $1$.

\begin{defi}
A matrix $v\in M_n(A)$ is called magic biunitary if all its rows and columns are partitions of the unity of $A$.
\end{defi}

A magic biunitary is indeed a biunitary, in the sense that both $v$ and its transpose $v^t$ are unitaries. The other word -- magic -- comes from a vague similarity with magic squares.

The basic example comes from the symmetric group $S_n$. Consider the sets of permutations $\{\sigma\in S_n\mid \sigma (j)=i \}$. When $i$ is fixed and $j$ varies, or vice versa, these sets form partitions of $S_n$. Thus their characteristic functions $v_{ij}\in \cc (S_n)$ form a magic biunitary.

Of particular interest is the ``universal'' magic biunitary matrix. This has coefficients in the universal algebra $\alg$ constructed by Wang in \cite{wa2}.

\begin{defi}
$\alg$ is the universal $\cc^*$-algebra generated by $n^2$ elements $u_{ij}$, subject to the magic biunitarity condition.
\end{defi}

In other words, we have the following universal property. For any magic biunitary matrix $v\in M_n(A)$ there is a morphism of $\cc^*$-algebras $\alg\to A$ mapping $u_{ij}\to v_{ij}$.

A more elaborate version of this property, to be discussed now, states that $A_{aut}(X_n)$ is a Hopf $\cc^*$-algebra, whose underlying quantum group is a kind of analogue of $S_n$.

The following definition is due to Woronowicz \cite{wo3}.

\begin{defi}
A unital Hopf $\cc^*$-algebra is a unital $\cc^*$-algebra $A$, together with a morphism of $\cc^*$-algebras $\Delta :A\to A\otimes A$, subject to the following conditions.

(i) Coassociativity condition: $(\Delta\otimes id)\Delta =(id\otimes\Delta )\Delta$.

(ii) Cocancellation condition: $\Delta (A)(1\otimes A)$ and $\Delta (A)(A\otimes 1)$ are dense in $A\otimes A$.
\end{defi}

The basic example is the algebra $\cc (G)$ of continuous functions on a compact group $G$, with $\Delta (\varphi): (g,h)\to\varphi (gh)$. Here coassociativity of $\Delta$ follows from associativity of the multiplication of $G$, and cocancellation in $\cc (G)$ follows from cancellation in $G$.

Another example is the group algebra $\cc^*(\Gamma )$ of a discrete group $\Gamma$. This is obtained from the usual group algebra $\cc [\Gamma]$ by a standard completion procedure. The comultiplication is defined on generators $g\in\Gamma$ by the formula $\Delta (g)=g\otimes g$.

In general, associated to a Hopf $\cc^*$-algebra $A$ are a compact quantum group $G$ and a discrete quantum group $\Gamma$, according to the heuristic formula $A=\cc (G)=\cc^*(\Gamma )$.

\begin{defi}
A coaction of $A$ on a finite set $X$ is a morphism of $\cc^*$-algebras $\alpha :\cc (X)\to \cc (X)\otimes A$, subject to the following conditions.

(i) Coassociativity condition: $(\alpha\otimes id)\alpha =(id\otimes\Delta )\alpha$.

(ii) Natural condition: $(\Sigma\otimes id)v=\Sigma(.)1$, where $\Sigma (\varphi)$ is the sum of values of $\varphi$.
\end{defi}

The basic example is with a group $G$ of permutations of $X$. Consider the action map $a:X\times G\to X$, given by $a(i,\sigma )=\sigma (i)$. The formula $\alpha\varphi =\varphi a$ defines a morphism of $\cc^*$-algebras $\alpha :\cc (X)\to\cc (X\times G)$. This can be regarded as a coaction of $\cc (G)$ on $X$.

In general, coactions of $A$ can be thought of as coming from actions of the underlying compact quantum group $G$. With this interpretation, the natural condition says that the action of $G$ must preserve the counting measure on $X$. This assumption cannot be dropped.

The following fundamental result is due to Wang \cite{wa2}.

\begin{theo}
(i) $\alg$ is a Hopf $\cc^*$-algebra, with comultiplication 
$\Delta (u_{ij})=\sum u_{ik}\otimes u_{kj}$.

(ii) The linear map $\alpha (\delta_j )=\sum\delta_i\otimes u_{ji}$ is a coaction of $\alg$ on $X_n=\{1,\ldots ,n\}$.

(iii) $\alg$ is the universal Hopf $\cc^*$-algebra coacting on $X_n$.
\end{theo}

The idea for proving (i) is that we can define $\Delta$ by using the universal property of $\alg$. Coassociativity is clear, and cocancellation follows from a result of Woronowicz in \cite{wo3}, stating that this is automatic whenever there is a counit and an antipode. But these can be defined by $\varepsilon (u_{ij})=\delta_{ij}$ and $S(u_{ij})=u_{ji}$, once again by using universality of $\alg$.

We know that the compact quantum group $G_n$ associated to $\alg$ is a kind of quantum analogue of the symmetric group $S_n$. In particular there should be an inclusion $S_n\subset G_n$. Here is the exact formulation of this observation, see Wang \cite{wa2} for details.

\begin{prop}
There is a Hopf $\cc^*$-algebra morphism $\pi_n :\alg\to \cc(S_n)$, mapping the generators $u_{ij}$ to the characteristic functions of the sets $\{ \sigma\in S_n\mid \sigma (j)=i\}$.
\end{prop}

The question is now whether $\pi_n$ is an isomorphism or not. For instance a $2\times 2$ magic biunitary must be of the following special form, where $p$ is a projection.
$$\begin{pmatrix}p&1-p\cr 1-p&p\end{pmatrix}$$

The algebra generated by $p$ is canonically isomorphic to $\cc^2$ if $p\neq 0,1$, and to $\cc$ if not. Thus the universal algebra $A_{aut}(X_2)$ is isomorphic to $\cc^2$, and $\pi_2$ is an isomorphism.

The map $\pi_3$ is an isomorphism as well, see \cite{gr} for a proof.

At $n=4$ we have the following example of magic biunitary matrix.
$$\begin{pmatrix}p&1-p&0&0\cr 1-p&p&0&0\cr 0&0&q&1-q\cr 0&0&1-q&q \end{pmatrix}$$

We can choose the projections $p,q$ such that the algebra $<p,q>$ they generate is infinite dimensional and not commutative. It follows that $A_{aut}(X_4)$ is infinite dimensional and not commutative as well, so $\pi_4$ cannot be an isomorphism.

\begin{prop}
For $n\geq 4$ the algebra $\alg$ is infinite dimensional and not commutative. In particular $\pi_n$ is not an isomorphism.
\end{prop}

This follows by gluing an identity matrix of size $n-4$ to the above $4\times 4$ matrix.

There is a quantum group interpretation here. Consider the compact and discrete quantum groups defined by the formula $A_{aut}(X_4)=\cc (G_4)=\cc^*(\Gamma_4)$. When $p,q$ are free the surjective morphism of $\cc^*$-algebras $A_{aut}(X_4)\to <p,q>$ can be thought of as coming from a surjective morphism of discrete quantum groups $\Gamma_4\to\z_2*\z_2$. This makes it clear that $\Gamma_4$ is infinite. Now $G_4$ being the Pontrjagin dual of $\Gamma_4$, it must be infinite as well.

See  Bichon \cite{bic}, Wang \cite{wa2}, \cite{wa4} and \cite{gr} for further speculations on this subject.

\section{Inner faithful representations}

We would like to find an explicit representation of $\alg$. As with any Hopf $\cc^*$-algebra, there is a problem here, because there are two notions of faithfulness.

Consider for instance a discrete subgroup $\Gamma$ of the unitary group $U(n)$. The inclusion $\Gamma\subset U(n)$ can be regarded as a unitary group representation $\Gamma\to U(n)$, and we get a $\cc^*$-algebra representation $\cc^*(\Gamma )\to M_n(\cc )$. This latter representation is far from being faithful: for instance its kernel is infinite dimensional, hence non-empty, when $\Gamma$ is an infinite group. However, the representation $\cc^*(\Gamma )\to M_n(\cc )$ must be ``inner faithful'' in some Hopf $\cc^*$-algebra sense, because the representation $\Gamma\to U(n)$ it comes from is faithful.

So, we are led to the following question. Let $H$ be a unital Hopf $\cc^*$-algebra, and let $\pi :H\to A$ be a morphism of $\cc^*$-algebras. If $\Gamma$ is the discrete quantum group associated to $H$ we know that $\pi$ corresponds to a unitary representation $\pi_i :\Gamma\to U(A)$. The question is: when is $\pi$ inner faithful, meaning that $\pi_i$ is faithful?

A simple answer is obtained by using the formalism of Kustermans and Vaes \cite{kv}. Associated to $H$ is a von Neumann algebra $H_{vN}$, obtained by a certain completion procedure. Now coefficients of $\pi$ belong to the dual algebra $\widehat{H}_{vN}$, and we can say that $\pi$ is inner faithful if these coefficients generate $\widehat{H}_{vN}$. This notion is used by Vaes in \cite{va}, and a version of it is used by Wang in \cite{wa2}.

In this paper we use an equivalent definition, from \cite{ver}.

\begin{defi}
Let $H$ be a unital Hopf $\cc^*$-algebra. A $\cc^*$-algebra representation $\pi :H\to A$ is called inner faithful if the $*$-algebra generated by its coefficients is dense in $H_{alg}^*$.
\end{defi}

Here $H_{alg}$ is the dense $*$-subalgebra of $H$ consisting of ``representative functions'' on the underlying compact quantum group, constructed by Woronowicz in \cite{wo3}. This is a Hopf $*$-algebra in the usual sense. Its dual complex vector space $H_{alg}^*$ is a $*$-algebra, with multiplication $\Delta^*$ and involution $*^*$. Finally, coefficients of $\pi$ are the linear forms $\varphi\pi$ with $\varphi\in A^*$, and the density assumption is with respect to the weak topology on $H_{alg}^*$. See e.g. the book of Abe \cite{abe} for Hopf algebras and \cite{ver} for details regarding this definition.

The main example is with a discrete group $\Gamma$. As expected, a representation $\cc^*(\Gamma )\to A$ is inner faithful if and only if the corresponding unitary group representation $\Gamma\to U(A)$ is faithful. Some other examples are discussed in \cite{ver}.

For how to use inner faithfulness see Vaes \cite{va}.

\begin{defi}
The character of a magic biunitary matrix $v\in M_n(A)$ is the sum of its diagonal entries $\chi (v)=v_{11}+v_{22}+\ldots +v_{nn}$.
\end{defi}

The terminology comes from the case where $v=u$ is the universal magic biunitary matrix, with coefficients in $A=\alg$. Indeed, the matrix $u$ is a corepresentation of $\alg$ in the sense of Woronowicz \cite{wo1}, and the element $\chi (u)$ is its character.

\begin{lemm}
Let $v\in M_n(A)$ be a magic biunitary matrix, with $n\geq 4$. Assume that there is a unital linear form $\varphi :A\to\cc$ such that
$$\varphi\left( \chi (v)^k\right)=\frac{1}{k+1}\begin{pmatrix}2k\cr k\end{pmatrix}$$
for any $k$. Then the representation $\pi:\alg\to A$ defined by $u_{ij}\to v_{ij}$ is inner faithful.
\end{lemm}

\begin{proof}
The numbers in the statement are the Catalan numbers, appearing as
multiplicities in representation theory of $SO(3)$. The result will follow
from the following fact from \cite{aut}. The finite dimensional irreducible corepresentations of $\alg$ can be arranged in a sequence $\{r_k\}$, such that their fusion rules are the same as those for representations of $SO(3)$.
$$r_k\otimes r_s=r_{|k-s|}+r_{|k-s|+1}+\ldots +r_{k+s}$$

Let $h:\alg\to\cc$ be the Haar functional, constructed by Woronowicz in \cite{wo1}. Consider also the character of the fundamental corepresentation of $\alg$.
$$\chi (u)=u_{11}+u_{22}+\ldots +u_{nn}$$

The Poincar\'e series of $\alg$ is defined by the following formula.
$$f(z)=\sum_{k=0}^\infty h\left( \chi (u)^k\right) z^k$$

By the above result, this is equal to the Poincar\'e series for $SO(3)$.
$$f(z)=\sum_{k=0}^\infty \frac{1}{k+1}\begin{pmatrix}2k\cr k\end{pmatrix}z^k$$

The assumption of the lemma says that the equality $\varphi\pi =h$ holds on all powers of $\chi (u)$. By linearity, this equality must hold on the algebra $<\chi (u)>$  generated by $\chi (u)$. Now by positivity of $h$ it follows that the restriction of $\pi$ to this algebra $<\chi (u)>$ is injective.

On the other hand, once again from fusion rules, we see that $\chi (u)$ generates the algebra of characters $\alg_{central}$ constructed by Woronowicz in \cite{wo1}.

Summing up, we know that $\pi$ is faithful on $\alg_{central}$.

Consider now the ``minimal model''  construction in \cite{ver}. This is the factorisation of $\pi$ into a Hopf $\cc^*$-algebra morphism $\alg\to H$, and an inner faithful representation $H\to A$.
$$\alg\to H\to A$$

Since $\pi$ is faithful on $\alg_{central}$, so is the map on the left. By Woronowicz's analogue of the Peter-Weyl theory in \cite{wo1} it follows that the map on the left is an isomorphism. Thus $\pi$ coincides with the map on the right, which is by definition inner faithful.
\end{proof}

It is possible to reformulate this result, by using notions from Voiculescu's free probability theory \cite{vdn}. A non-commutative $\cc^*$-probability space is a pair $(A,\varphi )$ consisting of a unital $\cc^*$-algebra $A$ together with a positive unital linear form $\varphi :A\to\cc$.

Associated to a self-adjoint element $x\in A$ is its spectral measure $\mu_x$. This is a probability measure on the spectrum of $x$, defined by the formula
$$\varphi (f(x ))=\int_\r f(t)\,d\mu_x(t).$$

This equality must hold for any continuous function $f$ on the spectrum of $x$. By density we can restrict attention to polynomials $f\in\cc [X]$, then by linearity it is enough to have this equality for monomials $f(t)=t^k$. We say that $\mu_x$ is uniquely determined by its moments,
$$\varphi (x^k)=\int_\r t^k\,d\mu_x(t).$$

The following notion plays a central role in free probability. See \cite{vdn},
page 26.

\begin{defi}
An element $x$ in a non-commutative $\cc^*$-probability space is called
semicircular if its spectral measure is
$d\mu_x(t)=(2\pi)^{-1}\sqrt{4-t^2}\,dt$ on $[-2,2]$, and $0$ elsewhere.
\end{defi}

In terms of moments, we must have the following equalities, for any $k$:
$$\varphi (x^k)=\frac{1}{2\pi}\int_{-2}^2t^k\sqrt{4-t^2}\,dt.$$

The integral is $0$ when $k$ is odd, and equal to a Catalan number when $k$ is even,
$$\varphi (x^{2k})= \frac{1}{k+1}\begin{pmatrix}2k\cr k\end{pmatrix}.$$

We get in this way a reformulation of the above lemma.

\begin{theo}
A magic biunitary matrix whose character has same spectral measure
as the square of a semicircular element produces an inner faithful
representation of Wang's algebra.
\end{theo}

The assumption $n\geq 4$ was removed, because it is superfluous. Indeed, for $n=1,2,3$ finite dimensionality of $\alg$ implies that the spectrum of any $\chi (v)$ is discrete.

\section{Geometric constructions}

Consider the Pauli matrices.
$$1=\begin{pmatrix}1&0\cr 0&1\end{pmatrix}\hskip 1cm
i=\begin{pmatrix}i&0\cr 0&-i\end{pmatrix}\hskip 1cm
j=\begin{pmatrix}0&1\cr -1&0\end{pmatrix}\hskip 1cm
k=\begin{pmatrix}0&i\cr i&0\end{pmatrix}$$

These satisfy the relations for quaternions $i^2=j^2=k^2=-1$, $ij=ji=-k$ etc.

To any $x\in SU(2)$ we associate the following matrix.
$$\begin{pmatrix}1\cr i\cr j\cr k\end{pmatrix}
x\begin{pmatrix}1&i&j&k\end{pmatrix}=
\begin{pmatrix}
x&xi&xj&xk\cr
ix&ixi&ixj&ixk\cr
jx&jxi&jxj&jxk\cr
kx&kxi&kxj&kxk\end{pmatrix}$$

Each row and each column of this matrix is an orthogonal basis of 
$M_2(\cc )\simeq\cc^4$ with respect to the inner product
\[\langle x,y\rangle=\frac12 Tr(xy^*),\]
since $i,j,k$ are sqew-adjoints.
Thus the matrix of corresponding orthogonal projections is a magic biunitary.

\begin{theo}
There is an inner faithful representation 
\[\pi :A_{aut}(X_4)\to \cc (SU(2), M_4(\cc ))\]
mapping the universal $4\times 4$ magic biunitary matrix to the $4\times 4$ matrix
$$v(x)=\begin{pmatrix}P_x&P_{xi}&P_{xj}&P_{xk}\cr
P_{ix}&P_{ixi}&P_{ixj}&P_{ixk}\cr
P_{jx}&P_{jxi}&P_{jxj}&P_{jxk}\cr
P_{kx}&P_{kxi}&P_{kxj}&P_{kxk}\end{pmatrix}$$
where for $y\in SU(2)$ we denote by $P_y$ the orthogonal projection onto the space $\cc y\subset M_2(\cc )$, and we regard it as a continuous function of $y$, with values in $M_2(M_2(\cc ))\simeq M_4(\cc )$. 
\end{theo}

\begin{proof}
We have to compute the character of $v=v(x)$.
$$\chi (v)=P_x+P_{ixi}+P_{jxj}+P_{kxk}$$

We make the convention that Greek letters designate quaternions in $\{
1,i,j,k\}$. We decompose $x$ as a sum with real coefficients $x=\sum
x_\alpha \alpha$. We have the following formula for $\chi (v)$.
$$\chi (v)=\sum_\alpha P_{\alpha x\alpha}$$

With the notations $\alpha\beta =(-1)^{N(\alpha ,\beta)}\beta\alpha$ and
$\alpha^2=(-1)^{N(\alpha )}$ we can compute $\alpha x\alpha$.
$$\alpha x\alpha =\sum_\beta (-1)^{N(\alpha ,\beta)+N(\alpha)}x_\beta\beta$$

Now using the above-mentioned canonical scalar product on $M_2(\cc)$, this gives the following formula for $P_{\alpha x\alpha}$, after cancelling
the $(-1)^{2N(\alpha)}=1$ term.
$$<P_{\alpha x\alpha}\beta ,\gamma >=(-1)^{N(\alpha ,\beta)+N(\alpha
  ,\gamma)}x_\beta x_\gamma$$

Now summing over $\alpha$ gives the formula of the character $\chi (v)$.
$$<\chi (v)\beta ,\gamma >=\sum_\alpha (-1)^{N(\alpha ,\beta )+N(\alpha
  ,\gamma )}x_\beta x_\gamma$$

The coefficient of $x_\beta x_\gamma$ can be computed by using the
multiplication table of quaternions.
$$\sum_\alpha (-1)^{N(\alpha ,\beta )+N(\alpha
  ,\gamma )}=4\delta_{\beta ,\gamma}$$

Thus $\chi (v)$ is a diagonal matrix, having the
numbers $4x_\beta^2$ on the diagonal.
$$\chi (v)={\mathrm{diag}} (4x_\beta^2)$$

Consider the linear form $\varphi =\int\otimes\, tr$, where the integral is
with respect to the Haar measure of $SU(2)$, and  $tr$ is the normalised trace of $4\times 4$ matrices, meaning $1/4$ times the usual trace. The moments of $\chi (v)$ with respect to $\varphi$ are computed as follows.
$$\int tr(\chi (v)^k)dx=4^{k-1}\sum_\beta\int x_\beta^{2k}dx$$

By symmetry reasons the four integrals are all equal, say to the first one.
$$\int tr(\chi (v)^k)dx=4^{k}\int x_1^{2k}dx$$

It follows that $\chi (v)$ has the same spectral measure as $4x_1^2$.
$$\mu_{\chi (v)}=\mu_{4x_1^2}$$

But the variable $2x_1$ is semicircular. This can be seen in many ways, for instance by direct computation, after identifying $SU(2)$
with the real sphere $S^3$, or by using the fact that $2x_1={\mathrm{Tr}}(x)$ is
the character of the fundamental representation of $SU(2)$, whose moments are
computed using Clebsch-Gordon rules. The result
follows now by applying theorem 2.1.
\end{proof}

The construction of $\pi$ has the following generalisation. Consider
the Clifford algebra $Cl(\r^s)$. This is a finite dimensional algebra, having a basis formed by products $e_{i_1}\ldots e_{i_k}$ with $1\leq
i_1<\ldots <i_k\leq s$, with multiplicative structure given by $e_i^2=-1$ and $e_ie_j=-e_je_i$ for $i\neq j$.

It is convenient to use the notation $e_I=e_{i_1}\ldots e_{i_k}$ with $I=(i_1,\ldots ,i_k)$. 

As an example, the Clifford algebra $Cl(\r^2)$ is spanned by the elements $e_\emptyset =1$, $e_1$, $e_2$ and $e_{12}=e_1e_2$. The generators $e_1,e_2$ are subject to the relations $e_1^2=e_2^2=-1$ and $e_1e_2=-e_2e_1$. Now these relations are satisfied by the Pauli matrices $i,j$, and the corresponding representation of $Cl(\r^2)$ turns to be faithful. That is, we have the following identifications.
$$e_\emptyset=\begin{pmatrix}1&0\cr 0&1\end{pmatrix}\hskip 1cm
e_1=\begin{pmatrix}i&0\cr 0&-i\end{pmatrix}\hskip 1cm
e_2=\begin{pmatrix}0&1\cr -1&0\end{pmatrix}\hskip 1cm
e_{12}=\begin{pmatrix}0&i\cr i&0\end{pmatrix}$$

We can label as well indices of $4\times 4$ matrices by elements of the set
$\{ \emptyset, 1,2,12\}$. With these notations, the representation in theorem
3.1 is given by $\pi(u_{IJ})=P_{e_Ixe_J}$.

The same formula works for an arbitrary number $s$.

\begin{theo}
There is a representation $\pi_n:A_{aut}(X_n)\to \cc (G_n,M_n(\cc
))$ mapping the universal $n\times n$ magic biunitary matrix to the $n\times
n$ matrix
$$v=\left( P_{e_Ixe_J}\right)_{IJ}$$
where $n=2^s$, the
unitary group of the Clifford algebra $Cl(\r^s)$ is denoted $G_n$, and the algebra of endomorphisms
of $Cl(\r^s)$ is identified with $M_{n}(\cc )$.
\end{theo}

The first part of proof of theorem 3.1 extends to this general situation. We
get that $\chi (v)$ is diagonal, with eigenvalues $\{
nx_I^2\}$. This doesn't seem to be related
to semicircular elements when $s\geq 3$. The representation $\pi_n$ probably comes
from an inner faithful representation of a quotient of $A_{aut}(X_n)$,
corresponding to a ``subgroup'' of the quantum permutation group.

\end{document}